\documentclass[12pt]{amsart}
\usepackage[T1]{fontenc}
\makeatletter

\newcommand{\LyX}{L\kern-.1667em\lower.25em\hbox{Y}\kern-.125emX\spacefactor1000}
\newcommand{\noun}[1]{\textsc{#1}}

\theoremstyle{plain}    
\newtheorem*{thm*}{Theorem} 
\theoremstyle{plain}    
\newtheorem{lem}{Lemma} 
\theoremstyle{remark}
\newtheorem*{rem*}{Remark}

\usepackage{thanks}
\usepackage{amssymb}
\makeatother

\begin{document}

\title{Representation of non periodic functions by trigonometric series with almost
integer frequencies}

\author{Gady Kozma and Alexander Olevski\v\i}

\address{School of Mathematical Sciences, Tel Aviv University, Ramat-Aviv, Israel 69978}

\address{e-mail address: gady@math.tau.ac.il, olevskii@math.tau.ac.il}

\maketitle
\begin{abstract}
Inspired by Men'shov's representation theorem, we prove that there exists a
sequence \( \left\{ \lambda (n)\right\} \subset \hat{\mathbf{R}} \), \( \lambda (n)=n+o(1) \),
\( n\in \mathbf{Z} \) such that any measurable (complex valued) function \( f \)
on \( \mathbf{R} \) can be represented as a sum of almost everywhere convergent
trigonometric series \( \sum c_{n}e^{i\lambda (n)x} \). 
\end{abstract}

\section{Introduction}

The Men'shov theorem (1940, see \cite[ch. XV]{bary}) states that every measurable
\( 2\pi  \)-periodic \( f \) can be represented as a sum
\[
f(x)=\sum _{n\in \mathbf{Z}}c_{n}e^{inx}\]
 of (nonunique) trigonometric series convergent a.e. on \( \mathbf{R} \). This
famous result has served as a starting point for many further investigations,
see \cite{TOyan} for a comprehensive survey and references. 

\Thanks{{} Partly supported by the Israel Science Foundation}

A non-periodic analog of the theorem is also known \cite{Davtjan}, where \( f \)
is expanded in a ``trigonometric integral'' involving all frequencies \( \lambda \in \hat{\mathbf{R}} \).

The aim of the present note is to show that by small perturbations of integers
one can get a universal spectrum of frequencies which allows to represent any
nonperiodic function on \( \mathbf{R} \) by pointwise convergent trigonometric
series.

\begin{thm*}
There exists a real sequence \( \Lambda =\left\{ \lambda (n)\right\} _{n=-\infty }^{\infty } \),

\begin{equation}
\label{lambdaisn}
\lambda (n)=n+o(1)
\end{equation}
 such that every measurable function \( f:\mathbf{R}\rightarrow \mathbf{C} \)
can be represented as a sum 
\begin{equation}
\label{fissum}
f(x)=\sum _{n\in \mathbf{Z}}c_{n}e^{i\lambda (n)x}
\end{equation}
convergent almost everywhere \\
Moreover, the perturbations \( \alpha (n)=\lambda (n)-n \) can be obtained
from an arbitrary pre-given sequence \( 0\neq \rho (k)=o(1) \) (\( k\in \mathbf{N} \))
by rearrangement with finite repetitions.
\end{thm*}
Convergence is understood with respect to symmetrical partial sums, i.e. \( \sum _{|\lambda (n)|<\eta }c_{n}e^{i\lambda (n)x} \),
\( \eta >0 \).

The idea of the proof below is to combine a recent result \cite{olevskii} about
approximation by polynomials with ``almost integer'' frequencies with the
elegant Körner version of Men'shov's technique, see \cite{Korner}.

\def \b {\bigskip}

\def \m {\medskip}

\def \s {\smallskip}

\b

{\bf Repr\'esentation de fonctions non-p\'eriodiques par des s\'eries trigonom\'etriques de fr\'equences presque enti\`eres.}

\b {\bf R\'esum\'e.} Inspir\'es par le th\'eor\`eme de repr\'esentation de Men'shov, nous prouvons qu'il existe une suite \( \left\{ \lambda (n)\right\} \subset \hat{\mathbf{R}} \), \( \lambda (n)=n+o(1) \), \( n\in \mathbf{Z} \) telle que toute fonction mesurable $f$ peut \^etre repr\'esent\'ee comme somme de s\'eries trigonom\'etriques convergentes presque partout.

\b {\bf Version francaise abr\'eg\'ee.} \s Le c\'el\`ebre th\'eor\`eme de Men'shov affirme que toute fonction p\'eriodique mesurable $f$ peut \^etre repr\'esent\'ee (de mani\`ere non-unique) comme somme de s\'eries trigonom\'etriques convergentes presque partout. \s Il y a de nombreuses versions et g\'en\'eralisations de ce r\'esultat (voir \cite{TOyan}). En particulier, un analogue non-periodique de ce th\'eor\`eme est connu (voir \cite{Davtjan}) dans lequel $f$ est decompos\'ee en "int\'egrale trigonom\'etrique" utilisant toutes les fr\'equences. \s Le but de cet article est de montrer que par des petites perturbations d'entiers on peut obtenir un spectre universel des fr\'equences qui permet de repr\'esenter toute fonction non-p\'eriodique sur \( \mathbf{R} \) par des s\'eries trigonom\'etriques convergentes point par point.

\m

{\bf Th\'eor\`eme.} {\it Il existe une suite r\'eelle \( \Lambda =\left\{ \lambda (n)\right\} _{n=-\infty }^{\infty } \), \( \lambda (n)=n+o(1) \) telle que toute fonction mesurable peut \^etre repr\'esent\'ee comme une somme \( f(x)=\sum _{n\in \mathbf{Z}}c_{n}e^{i\lambda (n)x} \) convergente p.p.

De plus, la perturbation peut \^etre obtenue \`a partir d'une suite arbitraire \( 0\neq \rho (k)=o(1) \) donn\'ee par r\'earrangement avec un nombre fini de r\'ep\'etitions. }

\s 

L'id\'ee de la preuve ci-dessous est de combiner un r\'esultat r\'ecent \cite{olevskii} sur l'approximation par des polyn\^omes de fr\'equences presque enti\`eres avec la version \'el\'egante de K\"orner \cite{Korner} des techniques de Men'shov.

\section{Preliminaries}

\subsection{Notations}

By a trigonometric polynomial we mean a finite sum 
\[
P(x)=\sum a_{k}e^{i\nu (k)x},\quad \cdots \nu (-1)<\nu (0)<\nu (1)<\cdots \subset \mathbf{R}\]
 The set \( \left\{ \nu (k)\right\}  \) (the spectrum of \( P \)) is denoted
by \( \mathrm{spec}\, P \). When the \( \nu (k) \) need to be integers, we
shall specify so explicitly, by saying that the polynomial has integer spectrum.
We will also use the following notations:
\begin{eqnarray*}
\deg P & = & \max |\nu (k)|\\
||P||_{A} & = & \sum |a_{k}|\\
P^{*}(x) & = & \sup _{\eta }\left| \sum _{|\nu (k)|<\eta }a_{k}e^{i\nu (k)x}\right| \\
||P||_{U} & = & \sup _{x\in \mathbf{R}}|P^{*}(x)|\\
P_{[N]}(x) & = & P(Nx)
\end{eqnarray*}
 \( \mathbf{m} \) will denote the Lebesgue measure on the line.

\subsection{External lemmas}

We use the following known results. The first is lemma 1 from \cite{Korner}

\begin{lem}
\label{korner}For every \( \delta >0 \) there exists some constant \( C(\delta ) \)
with satisfies the property that for every \( \epsilon >0 \) there exists a
trigonometric polynomial \( P_{\epsilon ,\delta } \) with integer spectrum
with the following properties:
\begin{enumerate}
\item \( \widehat{P_{\epsilon ,\delta }}(0)=0 \).
\item \( \left| \widehat{P_{\epsilon ,\delta }}(n)\right| <\epsilon \qquad \forall n \).
\item \( \left\Vert P_{\epsilon ,\delta }\right\Vert _{U}<C(\delta ) \).
\item \( \mathbf{m}\left( \left\{ x\in [0,2\pi ]\: :\: \left| P_{\epsilon ,\delta }(x)-1\right| \geq \epsilon \right\} \right) <\delta  \).
\end{enumerate}
\end{lem}
An inspection of Körner's proof of this lemma will show that it is possible
to choose these polynomials such that if \( \epsilon _{1}\leq \epsilon _{2} \)
and \( \delta _{1}\leq \delta _{2} \) then \( \deg P_{\epsilon _{1},\delta _{1}}\geq \deg P_{\epsilon _{2},\delta _{2}} \).
We shall denote \( d(\epsilon ,\delta ):=\deg P_{\epsilon ,\delta } \).

The next proposition is a slight variation of lemma 15 from \cite{Korner}.

\begin{lem}
\label{triv}Assume \( P \), \( Q \) are trigonometric polynomials, \( P \)
with integer spectrum, and \( N \) is a number satisfying \( N>2\deg Q \).
Then
\[
\left( P_{[N]}\cdot Q\right) ^{*}(x)\leq 2||\widehat{P}||_{\infty }\cdot ||Q||_{A}+\left| Q(x)\right| \cdot ||P||_{U}.\]

\end{lem}
The last result is taken from \cite{olevskii}:

\begin{lem}
\label{olevsk}If \( \sigma (n)=n+o(1) \), \( \sigma (n)\neq n \), then the
system \( \left\{ e^{i\sigma (n)x}\right\} _{n=-\infty }^{\infty } \) is complete
in \( L_{0}(\mathbf{R}) \) i.e. for every measurable function \( f:\mathbf{R}\rightarrow \mathbf{C} \)
there exists polynomials \( R_{k}=\sum c_{n,k}e^{i\sigma (n)x} \) such that
\( f(x)=\lim R_{k}(x) \) a.e. Further, this is true for any subsystem \( \left\{ e^{i\sigma (n)x}\right\} _{|n|>N} \).
\end{lem}

\section{Proof of Theorem}

\subsection{The sets \protect\( I_{l}\protect \) and the polynomials \protect\( R_{r,l}\protect \)}

Our first step is to define a fast-increasing sequence \( \eta _{l} \) such
that the sets \( I_{l}:=[-\eta _{l},-\eta _{l-1}]\cup [\eta _{l-1},\eta _{l}] \)
are ``large enough'' in the sense that the polynomials \( R_{k} \) in lemma
\ref{olevsk} can be constructed each with the spectrum supported on one \( I_{l} \).

More specifically, we shall construct polynomials \( \left\{ R_{r,l}\right\} _{0\leq l,\, |r|\leq l} \),
with spectrum in \( \left\{ \sigma (n)\right\} _{n=-\infty }^{\infty } \),
\( \sigma (n):=n+\rho (|n|) \) satisfying the following properties:

\begin{enumerate}
\item \( R_{r,l} \) have increasing spectra in the sense that if \( l_{1}>l_{2} \)
or \( l_{1}=l_{2} \) and \( r_{1}>r_{2} \) then \noun{
\[
\mathrm{spec}\, R_{r_{1},l_{1}}\subset \left\{ \xi \: :\: |\xi |>\deg R_{r_{2},l_{2}}\right\} \]
}
\item \( R_{r,l} \) approximates \( e^{i\sigma (r)x} \) in the sense 
\begin{equation}
\label{Rrl_property}
\mathbf{m}\left( \left\{ x\in [-l\pi ,l\pi ]\: :\: |R_{r,l}(x)-e^{i\sigma (r)x}|\geq \frac{1}{l^{2}}\right\} \right) <\frac{1}{l^{3}}
\end{equation}

\end{enumerate}
The self evident induction using lemma \ref{olevsk} will yield these polynomials.
We then define \( \eta _{l}:=\deg R_{l,l} \). Note that \( \mathrm{spec}\, R_{r,l}\subset I_{l} \).

\subsection{Construction of \protect\( \Lambda \protect \)}

We first need auxiliary sequences of numbers. \( \epsilon _{l} \), \( l\in \mathbf{N} \)
will be a sequence decreasing so fast that 
\begin{equation}
\label{eps_def}
\epsilon _{l}\cdot \max _{|r|\leq l}||R_{r,l}||_{A}<\frac{1}{l^{2}}
\end{equation}
With \( \epsilon _{l} \) we define \( d_{l}:=d(\epsilon _{l},l^{-3}) \) where
\( d(\epsilon ,\delta ) \) is defined after lemma \ref{korner}. Finally we
define \( b_{l} \) a sequence of integers satisfying \( b_{l}>b_{l-1}d_{l-1}+\eta _{l-1}+2\eta _{l} \).
If we now define sets \( I_{l,s}=I_{l}+sb_{l} \) then the sets \( I_{l,s} \)
will be disjoint for all \( l \) and \( 1\leq |s|\leq d_{l} \). We can now
define the sequence \( \Lambda  \) on the union \( J_{l}:=\bigcup _{1\leq |s|\leq d_{l}}I_{l,s} \)
as follows:
\[
\lambda (n+sb_{l})=\sigma (n)+sb_{l},\quad \sigma (n)\in I_{l},\; 1\leq |s|\leq d_{l}.\]
\( \Lambda  \) is thus defined on \( \left\{ n\: :\: \sigma (n)\in \bigcup _{l}J_{l}\right\}  \)
and clearly satisfies (\ref{lambdaisn}) and that \( \lambda (n)-n\in \left\{ \rho (k)\right\} _{k=0}^{\infty } \).
On the remaining \( n \)'s \( \Lambda  \) can be defined to be any arbitrary
sequence satisfying these two conditions --- for example \( \lambda (n)=\sigma (n) \).

\subsection{Representation of \protect\( f\protect \) -- definition of \protect\( c_{n}\protect \)}

Let now \( f:\mathbf{R}\rightarrow \mathbf{C} \) be any measureable function.
Our goal is to find coefficients \( c_{n} \) such that (\ref{fissum}) holds
a.e. We shall define successively both an increasing sequence \( \left\{ l(N)\right\}  \)
and blocks of coefficients corresponding to \( J_{l(N)} \) and set \( c_{n}=0 \)
if \( \sigma (n) \) does not belong to any \( J_{l(N)} \). Thus
\[
S_{N}:=\sum _{j=1}^{N}\sum _{\sigma (n)\in J_{l(N)}}c_{n}e^{i\lambda (n)x}\]
would be a subsequence of partial sums of the series (\ref{fissum}). Suppose
that \( N-1 \) steps are already done, so we have \( S_{N-1} \). We define
\( F_{N}:=f-S_{N-1} \).

\subsubsection{First approximation -- \protect\( G_{N}\protect \)}

We use lemma \ref{olevsk} to approximate \( F_{N} \) on the segment \( [-N\pi \, N\pi ] \)
with a uniform error of \( \delta _{N}=\frac{1}{NC((N+1)^{-3})} \) where the
function \( C(\delta ) \) is taken from lemma \ref{korner}; and a measure
error of \( \frac{1}{N^{2}} \) --- namely, find a polynomial \( G_{N}=\sum _{|r|<M_{N}}a_{r}e^{i\sigma (r)x} \)
satisfying 
\begin{equation}
\label{G-g}
\mathbf{m}\left( \left\{ x\in [-N\pi ,N\pi ]\: :\: |G_{N}(x)-F_{N}(x)|\geq \delta _{N}\right\} \right) <\frac{1}{N^{2}}
\end{equation}

\subsubsection{Second approximation -- \protect\( Q_{N}\protect \)}

Now we choose the integer \( l(N) \). We need it to be large enough, namely

\begin{equation}
\label{ln_def}
l(N)>N,\: \delta _{N}^{-1},\: M_{N},\: ||G_{N}||_{A},\: l(N-1)
\end{equation}
 \( Q_{N} \) is then defined as follows: 
\[
Q_{N}=\sum _{|r|<M_{N}}a_{r}R_{r,l(N)}.\]
The estimate of \( Q_{N}-G_{N} \) follows from (\ref{Rrl_property}):

\begin{eqnarray}
\mathbf{m}\left\{ x\in [-N\pi ,N\pi ]\: :\: |Q_{N}(x)-G_{N}(x)|\geq \delta _{N}\right\} \leq  &  & \nonumber \\
\sum _{|r|\leq M_{N}}\mathbf{m}\left\{ x\in [-N\pi \, N\pi ]\: :\: |R_{r,l(N)}(x)-e^{i\sigma (r)x}|\geq \frac{\delta _{N}}{||G_{N}||_{A}}\right\} \leq  &  & \nonumber \\
\sum _{|r|\leq M_{N}}\mathbf{m}\left\{ x\in [-l(N)\pi ,l(N)\pi ]\: :\: |R_{r,l(N)}(x)-e^{i\sigma (r)x}|\geq \frac{1}{l^{2}(N)}\right\} \leq  &  & \nonumber \\
\frac{2M_{N}+1}{l(N)^{3}}<\frac{3}{N^{2}} &  & \label{H-G} 
\end{eqnarray}
Notice also that
\[
\mathrm{spec}\, Q_{N}\subset I_{l(N)}\cap \left\{ \sigma (n)\right\} _{n=-\infty }^{\infty }\]
and finally that
\begin{equation}
\label{QnA_To_GnA}
||Q_{N}||_{A}\leq ||G_{N}||_{A}\cdot \max _{|r|\leq l(N)}||R_{r,l(N)}||_{A}.
\end{equation}

\subsubsection{Third approximation -- \protect\( H_{N}\protect \)}

The third approximation will come by multiplying \( Q_{N} \) with a certain
polynomial. We first use lemma \ref{korner}, with \( \delta =N^{-3} \) and
\( \epsilon =\epsilon _{l_{N}} \). and get a polynomial \( P_{N} \). Define

\[
H_{N}:=Q_{N}\cdot (P_{N})_{[b_{l(N)}]}\quad .\]
Notice that \( \deg P_{N}\leq d_{l(N)} \) and \( \widehat{P_{N}}(0)=0 \).
These two properties mean that \( \mathrm{spec}\, H_{N}\subset J_{l(N)} \).
Further, \( \mathrm{spec}\, Q_{N}\subset \left\{ \sigma (n)\right\} _{n=-\infty }^{\infty } \)
implies that \( \mathrm{spec}\, H_{N}\subset \Lambda  \). This allows us to
define 
\[
S_{N}:=S_{N-1}+H_{N}.\]

\subsection{Convergence of (\ref{fissum})}

We prove this in two stages: first that \( S_{N}(x)\rightarrow f(x) \) a.e.
and then that \( H_{N}^{*}(x)\rightarrow 0 \) a.e.

\subsubsection{}

Lemma \ref{korner} clause 4 gives

\[
\mathbf{m}\left( \left\{ x\in [-N\pi ,N\pi ]\: :\: \left| (P_{N})_{[b_{l(N)}]}(x)-1\right| \geq \epsilon \right\} \right) <N\delta =\frac{1}{N^{2}}.\]
Now, on the ``good'' set of \( x \) such that \( |P_{N}(b_{l(N)}x)-1|<\epsilon  \)
we can use (\ref{eps_def}) and get

\begin{eqnarray}
|H_{N}(x)-Q_{N}(x)| & \leq  & ||Q_{N}||_{A}\cdot |P_{N}(b_{l(N)}x)-1|\nonumber \\
 & \leq  & ||G_{N}||_{A}\cdot \max _{|r|\leq l(N)}||R_{r,l(N)}||_{A}\cdot \epsilon _{l(N)}\label{K-H} \\
 & \leq  & ||G_{N}||_{A}\cdot \frac{1}{l(N)^{2}}\leq \frac{1}{l(N)}\leq \delta _{N}.\nonumber 
\end{eqnarray}
Summing up (\ref{G-g}), (\ref{H-G}) and (\ref{K-H}) we get for \( f-S_{N}\equiv F_{N}-H_{N} \)

\begin{equation}
\label{f-S}
\mathbf{m}\left( \left\{ x\in [-N\pi \, N\pi ]\: :\: |f(x)-S_{N}(x)|\geq 3\delta _{N}\right\} \right) <\frac{5}{N^{2}}.
\end{equation}
which implies \( S_{N}(x)\rightarrow f(x) \) a.e. on \( \mathbf{R} \).

\subsubsection{}

Lemma \ref{triv} together with the estimate \( ||\hat{P}_{N}||_{\infty }<\epsilon _{l(N)} \)
gives us
\[
H_{N}^{*}(x)\leq 2\epsilon _{l(N)}\cdot ||Q_{N}||_{A}+\left| Q_{N}(x)\right| \cdot ||P_{N}||_{U}\qquad .\]
The first summand is \( <\frac{1}{N} \) by (\ref{eps_def}), (\ref{ln_def})
and (\ref{QnA_To_GnA}). Now, using (\ref{G-g}) and (\ref{H-G}) we get on
\( [-N\pi \, N\pi ] \) minus a set of measure \( \frac{4}{N^{2}} \) that 
\[
|Q_{N}(x)|\leq |F_{N}(x)|+2\delta _{N}\quad .\]
(\ref{f-S}) for \( N-1 \) implies \( |F_{N}(x)|\leq 3\delta _{N-1}=\frac{3}{(N-1)C(N^{-3})} \)
outside a set of measure \( \frac{5}{(N-1)^{2}} \) and finally \( ||P_{N}||_{U}<C(N^{-3}) \)
so the second summand is \( \leq \frac{5}{N-1} \) on \( [-N\pi \, N\pi ] \)
minus a set of measure \( <\frac{9}{(N-1)^{2}} \). It follows that that \( H^{*}_{N}(x)\rightarrow 0 \)
almost everywhere on \( \mathbf{R} \).\hfill $\square$

\begin{rem*}
It might be interesting to compare the approximation and expansion results.
The completeness theorem proved in \cite{olevskii} means that by \emph{arbitrary
small} perturbation of the integers one gets a spectrum \( \Lambda  \) which
is sufficient for approximation of any \( f\in L^{0}(\mathbf{R}) \) by linear
combination of \( e^{i\lambda x} \), \( \lambda \in \Lambda  \). In contrast
to that in the expansion theorem the perturbations can not decrease fast. In
particular, one can prove that in the theorem above it is impossible to construct
the sequence \( \Lambda  \) to satisfy the condition \( \lambda (n)=n+O(n^{-\epsilon }) \)
for some \( \epsilon >0 \). 
\end{rem*}

\end{document}